\documentclass{article}
\usepackage{latexsym,amsmath,amssymb,amscd}
\usepackage[mathscr]{eucal}
\usepackage[dvips]{color}
\begin{document}

\title{On Incompressible Vibrations of the Stratified Atmosphere on the Flat Earth }
\author{Tetu Makino \thanks{Professor Emeritus at Yamaguchi University, Japan; e-mail: makino@yamaguchi-u.ac.jp; Supported by JSPS KAKENHI Grant Number JP18K03371} }
\date{\today}
\maketitle

\newtheorem{Lemma}{Lemma}
\newtheorem{Proposition}{Proposition}
\newtheorem{Theorem}{Theorem}
\newtheorem{Definition}{Definition}
\newtheorem{Remark}{Remark}
\newtheorem{Corollary}{Corollary}
\newtheorem{Notation}{Notation}
\newtheorem{Assumption}{Assumption}

\numberwithin{equation}{section}

\begin{abstract}
We consider the vibrations and waves in the atmosphere under the gravitation on the flat earth. The stratified density distribution of the back ground equilibrium is supposed to touch the vacuum at the finite height of the stratosphere. We show  that incompressible motions are possible to create vibrations or progressive waves with countably many time frequencies of the vibrations or speeds of the wave propagation which accumulate to 0, say, of countably infinitely many slow and slow modes. 

 Eigenvalue problem  of Strum-Liouville type.


\end{abstract}

\section{Introduction}

We investigate the motion of the atmosphere on the flat earth under the constant gravitational force. We suppose that the motion is governed by the compressible Euler equations.

The pioneering mathematical investigation can be found as the paper `On the Vibrations of an Atmosphere' by Lord Rayleigh, 1890, \cite{Rayleigh}. He wrote

\begin{quote}
In order to introduce greater precision into our ideas respecting the behavior of the Earth's Atmosphere, it seems advisable to solve any problems that may present themselves, even though the search for simplicity may lead us to stray rather far from the actual question.
It is supposed here to consider the case of an atmosphere composed of gas which obeys Boyle's law, viz. such that the pressure is always proportional to the density. And in the first instance we shall neglect the curvature and rotation of the Earth, supposing that the strata of equal density are parallel planes perpendicular to the direction in which gravity acts.
\end{quote}

Our investigation in this article will be done under the same spirit as that of Lord Rayleigh quoted above, except for his starting point that the back ground state is supposed to be that of the isothermal stratified gas, say, $\rho=\rho_0\exp(-\mathsf{g}z/A)$, where $\rho$ is the density and $z$ is the height, $\rho_0, \mathsf{g}, A $ being positive constants. Instead we are interested in a back ground equilibrium like $\rho=C (z_+-z)^{\nu}$, where $z_+$ is the height of the stratosphere, $C, \nu (>1)$ being positive constants, which touches the vacuum $\rho=0$ 
on $z >z_+$ at the height $z=z_+$. We shall analyze the waves of the small amplitude around the back ground stratified equilibrium of this form assuming that the motion is done in the incompressible manner.  Therefore it might be better to call the problem as that of waves in a heterogeneous liquid, as \cite[Chapter IX, Section 235]{Lamb}, rather than the problem of waves in an atmosphere.  In contrast with the existing discussions, we should introduce the Lagrangian co-ordinate to treat the vacuum boundary as a free boundary.  Anyway the linearized analysis leads us to the investigation of the eigenvalue problem described by the so called Taylor-Goldstein equations. This equation is formally well known, e. g., as  \cite[p.378, (12)]{Lamb}, but, usually it was discussed in a bounded strap with both fixed boundaries at which the back ground density distribution behaves regularly and the usual Dirichlet boundary conditions are posed. If we consider it on the back ground which touches the vacuum at a finite height, it requires a careful analysis of the singular boundary
at the vacuum boundary in order to justify that the eigenvalue problem can be considered as that of the Strum-Liouville type.  
As the result we have vibrations with time frequency $\sqrt{\lambda_n}$ and progressive waves with speed $\sqrt{\lambda_n}/l$, $2\pi/l$ being the wave length in the horizontal direction, while
$\lambda_1 > \lambda_2 > \cdots > \lambda_n > \cdots \rightarrow 0$, in the linearized approximation level.
This analysis cannot be found in the existing literatures and is newly exhibited in this article.\\

Let us describe the problems more precisely. \\

We consider the motions of the atmosphere on the flat earth governed by the Euler equations
\begin{subequations}
\begin{align}
&\frac{\partial{\rho}}{\partial t}+\sum_{k=1}^3\frac{\partial}{\partial x^k}({\rho} v^k)=0, \label{1a} \\
&{\rho}\Big(\frac{\partial v^j}{\partial t}+\sum_{k=1}^3 v^k\frac{\partial v^j}{\partial x^k}\Big)
+\frac{\partial P}{\partial x^j}+{\rho}\frac{\partial\Phi}{\partial x^j}=0, \quad j=1,2,3,
 \label{1b} \\
&\Phi =\mathsf{g} x^3. \label{1c}
\end{align}
\end{subequations} 
Here $t\geq 0, \mathbf{x}=(x^1,x^2,x^3) \in \Omega:=\{\mathbf{x}\in\mathbb{R}^3 | 
x^3 >0\}$. The unknowns ${\rho} \geq 0, P $ are density and pressure, and
 $\displaystyle \mathbf{v}=v^1\frac{\partial}{\partial x^1}+v^2\frac{\partial}{\partial x^2}+v^3\frac{\partial}{\partial x^3}$ is the velocity field. $\mathsf{g}$ is a positive constant. The boundary condition is
\begin{equation}
{\rho} v^3=0 \quad\mbox{on}\quad x^3=0,
\end{equation}
and the initial condition is 
\begin{equation}
\rho=\overset{\circ}{\rho}(\mathbf{x}),\quad
\mathbf{v}=\overset{\circ}{\mathbf{v}}(\mathbf{x})=
(\overset{\circ}{v}^1(\mathbf{x}),
\overset{\circ}{v}^2(\mathbf{x}),
\overset{\circ}{v}^3(\mathbf{x}))\quad\mbox{at}\quad t=0.
\end{equation}\\

Moreover we suppose that the motion of the atmosphere is incompressible, that is, it holds that
\begin{equation}
\frac{D{\rho}}{Dt}:=\frac{\partial {\rho}}{\partial t}+\sum_{k=1}^3
v^k\frac{\partial{\rho}}{\partial x^k}=0, \label{3}
\end{equation}
which means that the density is constant along each stream line. Under the equation of continuity \eqref{1a}, this is equivalent to
\begin{equation}
{\rho} \sum_{k=1}^3\frac{\partial v^k}{\partial x^k}=0. \label{4}
\end{equation}\\

Suppose that there is fixed a flat equilibrium ${\rho}=\bar{\rho}, P=\bar{P}, \mathbf{v}=\mathbf{0}$ which satisfy \eqref{1a}\eqref{1b} with \eqref{1c}, where $\bar{\rho}, \bar{P}$ are functions of $z=x^3$ only. Namely we suppose that
\begin{equation}
\frac{d\bar{P}}{dz}+\mathsf{g}\bar{\rho}=0\qquad\mbox{for}\qquad z >0.
\end{equation}
We assume that $\bar{\rho}\in C^1([0, +\infty[)$ and
$\{\bar{\rho} >0\}=\{ 0\leq z<z_+\}$, where $0<z_+\leq +\infty$. We suppose
\begin{equation}
\frac{d\bar{\rho}}{dz} <0\qquad\mbox{for}\qquad 0\leq z <z_+.
\end{equation}\\

We consider the perturbation $\xi^j=\delta x^j, j=1,2,3, \delta{\rho},
\delta P$ at this fixed equilibrium. 

\begin{Notation}
We use the Lagrangian coordinate which will be denoted by the diversion of the letter $\mathbf{x}$ of the Eulerian coordinate. Instead, hereafter, we shall use the symbols
$\underline{\mathbf{x}}=(\underline{x}^1,\underline{x}^2, \underline{x}^3)$ to denote the original Eulerian co-ordinates so that $\underline{\mathbf{x}}=
\mathbf{x}+\mbox{\boldmath$\xi$}(t, \mathbf{x})$, $\mbox{\boldmath$\xi$}$ being
$(\xi^1,\xi^2,\xi^3)$.
\end{Notation}

Of course we put
\begin{equation}
\xi^1=\xi^2=\xi^3=0\qquad\mbox{at}\qquad t=0
\end{equation}
The boundary condition is
\begin{equation}
\xi^3=0\qquad\mbox{on}\qquad x^3=0,
\end{equation}
and the initial condition is
\begin{equation}
\frac{\partial \xi^j}{\partial t}\Big|_{t=0}=\overset{\circ}{v}^j(\mathbf{x}),\qquad j=1,2,3.
\end{equation}

Here let us recall the definition of the Euler perturbation $\delta Q$ and the Lagrange perturbation 
$\Delta Q$ of a quantity $Q$:
\begin{align*}
&\Delta Q(t,\mathbf{x})=Q(t, \mbox{\boldmath$\varphi$}(t,\mathbf{x}))-\bar{Q}(\mathbf{x}), \\
&\delta Q(t,\mathbf{x})=Q(t,\mbox{\boldmath$\varphi$}(t,\mathbf{x}))-\bar{Q}
(\mbox{\boldmath$\varphi$}(t,\mathbf{x})),
\end{align*}
where $\underline{\mathbf{x}}=\mbox{\boldmath$\varphi$}(t,\mathbf{x})=\mathbf{x}+\mbox{\boldmath$\xi$}(t,\mathbf{x})$ is the steam line given by 
$$\frac{\partial}{\partial t}\mbox{\boldmath$\varphi$}(t,\mathbf{x})=\mathbf{v}(t,\mbox{\boldmath$\varphi$}(t,\mathbf{x})),
\quad \mbox{\boldmath$\varphi$}(0,\mathbf{x})=\mathbf{x}. $$

We consider that the Lagrangian co-ordinates variables $\mathbf{x}=(x^1,x^2,x^3)$
runs on the fixed domain
\begin{equation}
\Pi:=\{ \mathbf{x}\in\mathbb{R}^3 \quad |\quad  \overset{\circ}{\rho}(\mathbf{x}) >0 \}.
\end{equation}
Note that
\begin{equation}
\overset{\circ}{\rho}(\mathbf{x})=\bar{\rho}(\mathbf{x})+\Delta\rho|_{t=0}(\mathbf{x}) =\bar{\rho}(\mathbf{x})+\delta\rho|_{t=0}(\mathbf{x}).
\end{equation}\\

Then the equation \eqref{1b} reads
\begin{equation}
\frac{\partial^2\mbox{\boldmath$\xi$}}{\partial t^2}+
\mathsf{g}\frac{\delta{\rho}}{\bar{\rho}+\delta{\rho}}\mathbf{e}_3+
\frac{1}{\bar{\rho}+\delta{\rho}}J^{-1}\nabla\delta P=0, \label{NLEq}
\end{equation}
where 
$$ \mbox{\boldmath$\xi$}=\sum_{k=1}^3\xi^k\frac{\partial}{\partial x^k},\qquad
\nabla \delta P=\sum_{k=1}^3\frac{\partial \delta P}{\partial x^k}\frac{\partial}{\partial x^k},
\qquad 
\mathbf{e}_3=\frac{\partial}{\partial x^3}, $$
and $J^{-1}=(({J}^{-1})_j^k)_{k,j}=(\partial x^k/\partial \underline{x}^j)_{k,j}$
 is the inverse matrix of 
\begin{equation}
J=(J_k^j)_{j,k}=\Big(\frac{\partial \underline{x}^j}{\partial x^k}\Big)_{j,k}=
\Big(\delta_k^j+\frac{\partial\xi^j}{\partial x^k}\Big)_{j,k}.
\end{equation}
Now the assumption \eqref{4}
$$0=\sum_k\frac{\partial v^k}{\partial\underline{x}^k}
=\sum_{k,j}(J^{-1})_k^j\frac{\partial}{\partial t}\frac{\partial\xi^k}{\partial x^j}=
\sum_{j,k}(J^{-1})_k^j\frac{\partial}{\partial t}J_j^k
$$
says that
$$
\mathrm{tr}\Big(\frac{\partial J^{-1}}{\partial t} J\Big)=0,
$$
and, equivalently, that
\begin{equation}
\mathrm{det}J^{-1}=\frac{1}{\mathrm{det}J}=1 \label{1.15}
\end{equation}
holds for $t \geq 0$ and $\mathbf{x} \in \Pi$, since
$J|_{t=0}=I$, $I=(\delta_k^j)_{j,k}$ being the unit matrix.\\

Note that we consider the equation \eqref{NLEq} on the fixed domain
$[0,T[\times\Pi $
and the dominator $\bar{\rho}+\delta{\rho}$ in the terms of the equation \eqref{NLEq} should read 
\begin{equation}
\bar{\rho}+\delta{\rho}=\bar{\rho}(\mathbf{x}+\mbox{\boldmath$\xi$}(t,\mathbf{x}))+\delta{\rho}(t,\mathbf{x}).\label{0111}
\end{equation}
Moreover we see that the Lagrangian perturbation $\Delta{\rho}$ is independent of $t$ thanks to \eqref{3}. Therefore by the definition we have
\begin{align}
\delta{\rho}(t,\mathbf{x})&=-(\bar{\rho}(\mathbf{x}+\mbox{\boldmath$\xi$}(t,\mathbf{x}))-\bar{\rho}(\mathbf{x})) +\Delta{\rho}|_{t=0}(\mathbf{x})  \nonumber \\
&=-(\bar{\rho}(\mathbf{x}+\mbox{\boldmath$\xi$}(t,\mathbf{x}))-\bar{\rho}(\mathbf{x})) +\delta{\rho}|_{t=0}(\mathbf{x}). \label{deltarho}
\end{align}
Accordingly we should read \eqref{0111} as
\begin{equation}
\bar{\rho}+\delta{\rho}=\bar{\rho}(\mathbf{x})+\delta\rho|_{t=0}(\mathbf{x})
=\overset{\circ}{\rho}(\mathbf{x}),
\end{equation}
and the dominator $\bar{\rho}+\delta{\rho}$ in the terms of the equation
\eqref{NLEq} does not contain unknown functions and is nothing but the initial density distribution. \\

If we insert $P=A\bar{\rho}^{\gamma}$ with $1<\gamma <2$ into $dP/dz=-\mathsf{g}\bar{\rho}$,
then we get
$$
\bar{\rho}(z)=C((z_+-z)\vee 0)^{\nu}=
\begin{cases}
C(z_+-z)^{\nu} \qquad &(0\leq z <z_+) \\
0 \qquad &(z_+\leq z)
\end{cases}
$$
Here $z_+$ is an arbitrary finite positive number, $\displaystyle \nu:=\frac{1}{\gamma-1}$,  and 
$$C=\Big(\frac{(\gamma-1)\mathsf{g}}{A\gamma}\Big)^{\nu}. $$ 
Keeping in mind this case, we consider an equilibrium ${\rho}=\bar{\rho}(z)$ with $P=\bar{P}(z)$ which satisfies the following conditions:

1) $\{\bar{\rho} >0\}=\{0 \leq z <z_+\}$ with a finite positive $z_+$, and
\begin{equation}
\bar{\rho} \in C^1([0, +\infty[)\cap C^{\infty}([0, z_+[);
\end{equation}

2) It holds
\begin{equation}
\frac{d\bar{\rho}}{dz} <0 \qquad\mbox{for}\qquad 0\leq z <z_+;
\end{equation}

3) $\bar{\rho}$ is analytic at $z=0$ and 
\begin{equation}
\bar{\rho} = C_{\bar{\rho}}(z_+-z)^{\nu}(1+\Lambda(z_+-z) )
\end{equation}
as $z \rightarrow z_+-0$, where $C_{\bar{\rho}}$ is a positive constant and $\Lambda$ is an analytic function near $0$ such that $\Lambda(0)=0$.

Hereafter we shall use the following notation so that $\Lambda(z_+-z)=[z_+-z]_1$:

\begin{Notation}
Let $K$ be a non-negative integer. Then the symbol $[X]_K$ stands for  various convergent power series of the form
$\sum_{k\geq K}a_kX^k$.
\end{Notation}

Note that 
\begin{align*}
\frac{d\bar{P}}{dz}&=-\mathsf{g}\bar{\rho}=\mathsf{g}C_{\bar{\rho}}(z_+-z)^{\nu}(1+[z_+-z]_1), \\
\frac{d\bar{\rho}}{dz}&=-\nu C_{\bar{\rho}}(z_+-z)^{\nu-1}(1+[z_+-z]_1),
\end{align*}
and therefore
$$\frac{d\bar{P}}{d\bar{\rho}}=\frac{\mathsf{g}}{\nu}(z_+-z)(1+[z_+-z]_1)
$$
so that 
$$-\infty < \frac{d}{dz}\Big(\frac{d\bar{P}}{d\bar{\rho}}\Big)\Big|_{z=z_+-0}=-\frac{\mathsf{g}}{\nu} <0,
$$
that is, the boundary $z=z_+$ is a so called `physical vacuum boundary'. See the review \cite{JJM}. Here $\displaystyle \frac{d\bar{P}}{d\bar{\rho}} $
means 
$\displaystyle \Big(\frac{d\bar{\rho}}{dz}\Big)^{-1}\Big(\frac{d\bar{P}}{dz}\Big)=\mathsf{g}H[\bar{\rho}], H[\bar{\rho}]=\Big(-\frac{d}{dz}\log\bar{\rho}\Big)^{-1}$ being the so called `density scale height'.\\

In this article we fix such an equilibrium $\bar{\rho}$ with $\bar{P}$ and investigate small perturbations around this fixed equilibrium.

\section{Linearized problem}

\subsection{}

Let us derive the linearized approximation of the equations.

The linearized approximation of \eqref{deltarho} reads
\begin{equation}
\delta{\rho}=-(\nabla|\bar{\rho}\mbox{\boldmath$\xi$}) 
+\delta{\rho}|_{t=0}. \label{8}
\end{equation}
 Here we have used the fact that \eqref{4} implies 
\begin{equation}
(\nabla|\mbox{\boldmath$\xi$})=0 \label{9}
\end{equation}
in the linearized approximation. The linearized approximation of the equation
\eqref{NLEq} is
\begin{equation}
\frac{\partial^2\mbox{\boldmath$\xi$}}{\partial t^2}+\mathsf{g}\frac{\delta{\rho}}{\bar{\rho}}\mathbf{e}_3+
\frac{1}{\bar{\rho}}\nabla\delta P=0. \label{10}
\end{equation}

\subsection{}

Let us consider particular solutions of \eqref{8}\eqref{9}\eqref{10} 
of the following two types of the form:\\

{\bf Type (1):}
\begin{subequations}
\begin{align}
\mbox{\boldmath$\xi$}&=
\begin{bmatrix}
u(z)\cos lx \\
0 \\
w(z)\sin lx
\end{bmatrix}
\sin \sqrt{\lambda}t, \label{2.1a}\\
\delta P&=\delta\check{P}(z)\sin lx\sin\sqrt{\lambda}t. \label{2.1b}
\end{align}
\end{subequations}
Here $x=x^1, z=x^3$, $\lambda$ is a positive constant and $l$ is a non-negative constant.
In this type, taking $\delta{\rho}|_{t=0}=0$, that is, $\overset{\circ}{\rho}=\bar{\rho}$, we put
\begin{equation}
\delta{\rho}=-(\nabla|\bar{\rho}\mbox{\boldmath$\xi$}). 
\end{equation}\\

{\bf Type (2):}
\begin{subequations}
\begin{align}
&\mbox{\boldmath$\xi$}=
\begin{bmatrix}
u(z)(\cos(lx-\sqrt{\lambda}t)-\cos lx) \\
0 \\
w(z)(\sin(lx-\sqrt{\lambda}t)-\sin lx)
\end{bmatrix}, \label{Type1a}\\
&\delta P=\delta\check{P}(z)(\sin(lx-\sqrt{\lambda}t)-\sin lx) \label{Type1b}
\end{align}
\end{subequations}
instead of \eqref{2.1a}\eqref{2.1b}. In this type we should specify the initial perturbation of the density
$\Delta{\rho}|_{t=0}=\delta{\rho}|_{t=0}$ by
\begin{equation}
\delta{\rho}|_{t=0}=\Big[l\bar{\rho} u -\frac{d}{dz}(\bar{\rho} w)\Big]\sin lx. \label{0207}
\end{equation}
Then we have
\begin{align}
\delta{\rho}&=-(\nabla|\bar{\rho} \mbox{\boldmath$\xi$})+\delta{\rho}|_{t=0} \nonumber \\
&=\Big[l\bar{\rho} u -\frac{d}{dz}(\bar{\rho} w)\Big]\sin(lx -\sqrt{\lambda}t), \label{0208}
\end{align}\\

 Then, for each type, the equations \eqref{9} and \eqref{10} are reduced to
\begin{equation}
-lu+\frac{dw}{dz}=0, \label{14c}
\end{equation}
and
\begin{subequations}
\begin{align}
-\lambda u +\frac{l}{\bar{\rho}}\delta\check{P}&=0, \label{14a}\\
-\lambda w-\mathsf{g}\frac{1}{\bar{\rho}}\frac{d\bar{\rho}}{dz}w+\frac{1}{\bar{\rho}}\frac{d}{dz}\delta\check{P}&=0. \label{14b}
\end{align}
\end{subequations}\\

Note that $\Pi=\{ 0 \leq z <z_+\}$ for the Type (1), but
\begin{align*}
\Pi&=\{ \bar{\rho}(z)-
\frac{d\bar{\rho}}{dz}w(z)\sin lx >0 \} \\
&=\{ 0\leq z <z_+\quad\mbox{and}\quad
1-\frac{1}{\bar{\rho}}\frac{d\bar{\rho}}{dz}w(z)\sin lx > 0\}
\end{align*}
for the Type (2), since \eqref{14c} implies that 
\eqref{0208} reads 
\begin{equation}
\delta\rho=-\frac{d\bar{\rho}}{dz}w(z)\sin(lx -\sqrt{\lambda}t).
\end{equation}\\

The set of equations \eqref{14a}\eqref{14c} is equivalent to 
\begin{equation}
\delta\check{P}=\frac{\lambda}{l}\bar{\rho} u=\frac{\lambda}{l^2}\bar{\rho}\frac{dw}{dz},
\end{equation}
and therefore \eqref{14b} turns out to be
\begin{equation}
\frac{d}{dz}\Big(\bar{\rho}\frac{dw}{dz}\Big)+
\frac{l^2}{\lambda}\Big(-\mathsf{g}\frac{d\bar{\rho}}{dz}-\bar{\rho}\lambda\Big)w=0. \label{16}
\end{equation}
Putting
\begin{equation}
\mathcal{N}:=\sqrt{-\frac{\mathsf{g}}{\bar{\rho}}\frac{d\bar{\rho}}{dz}},
\end{equation}
 we can write \eqref{16} as
\begin{equation}
\frac{d}{dz}\Big(\bar{\rho}\frac{dw}{dz}\Big)+\frac{l^2}{\lambda}\bar{\rho}
(\mathcal{N}^2-\lambda)w=0, \label{TG}
\end{equation}
which is called the  `Taylor-Goldstein equation'. The boundary condition is
\begin{equation}
w=0\qquad\mbox{at}\qquad z=0 \label{TGBC}
\end{equation}\\

The equation \eqref{TG} can be written as 
\begin{equation}
-\frac{d}{dz}\Big(\bar{\rho}\frac{dw}{dz}\Big)+l^2\bar{\rho} w=\frac{1}{\lambda}\mu(z)w, \label{0211}
\end{equation}
where
\begin{equation}
\mu(z):=-\mathsf{g}l^2\frac{d\bar{\rho}}{dz} \label{0212}
\end{equation}
Note that $\mu(z) >0$ for $0\leq z <z_+$. Let us perform the Liouville transformation on the eigenvalue problem
\eqref{0211}, where the eigenvalue is $\displaystyle \frac{1}{\lambda}$. Namely, putting
\begin{subequations}
\begin{align}
\zeta&=\int_0^z\sqrt{\frac{\mu}{\bar{\rho}}}(z')dz', \\
\upsilon&=(\bar{\rho}\mu)^{\frac{1}{4}}w, \\
q&=\frac{\bar{\rho}}{\mu}\Big[ l^2+\frac{1}{4}\frac{d^2}{dz^2}\log\Big(-\bar{\rho}\frac{d\bar{\rho}}{dz}\Big)+ \nonumber \\
&-\frac{1}{16}\Big(\frac{d}{dz}\log\Big(-\bar{\rho}\frac{d\bar{\rho}}{dz}\Big)\Big)^2+
\frac{1}{4}\Big(\frac{d}{dz}\log \bar{\rho}\Big)
\Big(\frac{d}{dz}\log\Big(-\bar{\rho}\frac{d\bar{\rho}}{dz}\Big)\Big)\Big],
\end{align}
\end{subequations}
we transform \eqref{0211} to the standard form
\begin{equation}
-\frac{d^2\upsilon}{d\zeta^2}+q\upsilon=\frac{1}{\lambda}\upsilon. \label{0214}
\end{equation}
By the property 3) of $\bar{\rho}$ we have, as $z \rightarrow z_+-0$, 
\begin{equation}
\frac{\mu}{\bar{\rho}}=-\nu\mathsf{g}l^2(z_+-z)^{-1}(1+[z_+-z]_1), \label{0215}
\end{equation}
therefore 
\begin{equation}
\zeta_+:=\int_0^{z_+}\sqrt{\frac{\mu}{\bar{\rho}}}(z)dz
\end{equation}
is finite and
\begin{equation}
\zeta_+-\zeta=2\sqrt{\nu\mathsf{g}l}\sqrt{z_+-z}(1+[z_+-z]_1). \label{0217}
\end{equation}
Thus the $z$-interval $[0, z_+]$ is mapped onto the $\zeta$-interval $[0, \zeta_+]$.

On the other hand, be a tedious calculation, we have
$$q=\frac{(2\nu-1)(2\nu-3)}{16\nu\mathsf{g}l^2}
(z_+-z)^{-1}(1+[z_+-z]_1) $$
if $\nu \not=3/2$ and 
$$q=[z_+-z]_0 $$
if $\nu=3/2$, or , thanks to \eqref{0217},
\begin{equation}
q=\frac{(2\nu-1)(2\nu-3)}{4}(\zeta_+-\zeta)^{-2}
(1+[(\zeta_+-\zeta)^2]_1) \label{0218}
\end{equation}
if $\nu \not=3/2$ and
\begin{equation}
q=[(\zeta_+-\zeta)^2]_0 \label{0219}
\end{equation}
if $\nu=3/2$. Note that $\displaystyle \frac{(2\nu-1)(2\nu-3)}{4} > -\frac{1}{4}$ for $\nu >1$. Since $q \in C([0, \zeta_+[)$, we can claim the following

\begin{Proposition}
There are constants 
$\displaystyle K_0 >-\infty, K_1 >-\frac{1}{4}$ such that
\begin{equation}
q(\zeta) \geq K_0+\frac{K_1}{(\zeta_+-\zeta)^2}\qquad\mbox{on}\qquad
0\leq \zeta <\zeta_+. \label{0220}
\end{equation}
\end{Proposition}

Consequently we can claim the following

\begin{Theorem}\label{Th.1}
The operator $\displaystyle -\frac{d^2}{d\zeta^2}+q$ defined on $C_0^{\infty}(]0,\zeta_+[)$
$[\!($ or The operator $\displaystyle -\frac{d}{dz}\bar{\rho}\frac{d}{dz}+l^2\bar{\rho}$ defined on $C_0^{\infty}(]0, z_+[)$ $)\!]$
 admits the Friedrichs extension, which is a self-adjoint operator, in the Hilbert space
$L^2([0,\zeta_+])$ $[\!($ in the Hilbert sapce $\displaystyle  L^2([0,z_+], \mathsf{g}l^2\Big(-\frac{d\bar{\rho}}{dz}\Big)dz), respectively )\!] $. The spectrum of the self-adjoint operator consists of eigenvalues  
$$0<\frac{1}{\lambda_1} <\frac{1}{\lambda_2 }<\cdots <\frac{1}{\lambda_n} <\cdots \rightarrow +\infty.$$
\end{Theorem}
For a proof see, e.g., \cite[Kapiter VII]{CourantH}. Let us fix an eigenvalue $\frac{1}{\lambda_n}$ and an associated eigenfunction $\upsilon=\upsilon_n(\zeta)$. We are going to prove the following
\begin{Proposition}\label{Prop.1}
There is a constant $a$ such that $a\not=0$ and
\begin{equation}
\upsilon_n(\zeta)=a(\zeta_+-\zeta)^{\frac{2\nu-1}{2}}(1+[(\zeta_+-\zeta)^2]_1) \label{0221}
\end{equation}
\end{Proposition}

In order to prove the Proposition \ref{Prop.1}, we shall use the following
\begin{Lemma}\label{Lem.1}
Let us consider the equation
\begin{equation}
-\frac{d^2y}{dx^2}+V(x)y=0,\label{0222}
\end{equation}
where
\begin{equation}
V(x)=\frac{1}{x^2}(K+[x^{\beta}]_1). \label{0223}
\end{equation}
with constants $K, \beta$ such that $4K+1>0, \beta >0$. Put
\begin{equation}
\alpha_{\pm}:=\frac{1}{2}(1\pm\sqrt{4K+1}).\label{0224}
\end{equation}
Then there is a fundamental system of solutions $y=\varphi_1(x), y=\varphi_2(x)$ of the equation \eqref{0222} such that
\begin{subequations}
\begin{align}
\varphi_1(x)&=x^{\alpha_+}(1+[x^{\beta}]_1), \\
\varphi_2(x)&=x^{\alpha_-}(1+[x^{\beta}]_1)+h(\log x)\varphi_1(x).
\end{align}
\end{subequations}
Here $h$ is a constant, which may $=0$ or $\not=0$ generally, but $h=0$ if $\alpha_+-\alpha_-$ is not an integer.
\end{Lemma}

For a proof, see e.g., \cite[Chapter 1, \S 3]{Yosida}, \cite[Chapter 4]{CoddingtonL}. 
Actually, when $\beta \not=1$, we can reduce the proof to that of the case with $\beta=1$, say, the case of the regular singular point, by the change of variables $x \leadsto x_{\natural}:=x^{\beta}, y \leadsto
y_{\natural}:=x^{\frac{\beta-1}{2}}y$.\\

Proof of Proposition \ref{Prop.1}. Let us suppose first that $\nu\not=3/2$. When
$K=(2\nu-1)(2\nu-3)/4$, then $\alpha_{\pm}$ of Lemma \ref{Lem.1} turn out to be
$$\alpha_+=\frac{1}{2}(2\nu-1),\qquad \alpha_-=\frac{1}{2}(-2\nu+3). $$
So, by the Lemma \ref{Lem.1} we have a fundamental system of solutions $\upsilon^{(1)}, \upsilon^{(2)}$ of \eqref{0214} such that
\begin{align*}
\upsilon^{(1)}(\zeta)&=(\zeta_+-\zeta)^{\frac{2\nu-1}{2}}(1+[(\zeta_+-\zeta)^2]_1), \\
\upsilon^{(2)}(\zeta)&=(\zeta_+-\zeta)^{\frac{-2\nu+3}{2}}(1+[(\zeta_+-\zeta)^2]_1)+ \\
&+h(\log (\zeta_+-\zeta))\upsilon^{(1)}(\zeta).
\end{align*}
But $\upsilon^{(2)}$ does not belong to the domain of the Friedrichs self-adjoint extension of $-\frac{d^2}{d\zeta^2}+q \upharpoonright C_0^{\infty} $ in $L^2$, since
$d\upsilon^{(2)}/d\zeta \sim \frac{2\nu-3}{2}(\zeta_+-\zeta)^{\frac{-2\nu+1}{2}}
\not\in L^2$. Thus $\upsilon_n=a\upsilon^{(1)}$ with $a\not=0$. Next we suppose $\nu=3/2$. Then the boundary point $z=z_+$ is a regular boundary point, if we reduce the equation \eqref{0214} in the variable $(\zeta_+-\zeta)^2$. Therefore there is a fundamental system of solutions $\upsilon^{(1)}, \upsilon^{(2)}$ such that
\begin{align*}
\upsilon^{(1)}(\zeta)&=(\zeta_+-\zeta)(1+[(\zeta_+-\zeta)^2]_1, \\
\upsilon^{(2)}(\zeta)&=1+[(\zeta_+-\zeta)^2]_1.
\end{align*}
Since a function $\upsilon$ belongs to the domain of the Firedrichs extension if and only if it satisfies the Dirichret boundary conditions $\upsilon=0$ both at $\zeta=0$ and $\zeta=\zeta_+$, it should hold that $\upsilon_n=a\upsilon^{(1)}$ with $a\not=0$. This completes the proof of Proposition \ref{Prop.1} $\square$ \\

By the argument in the proof of Proposition \ref{Prop.1} we can claim the following

\begin{Corollary}
The eigenvalues $1/\lambda_n$ of Theorem \ref{Th.1} are simple.
\end{Corollary}

Note that, if and only if the continuation of the solution $\upsilon^{(1)}$ to the left hits $0$ at the regular boundary point $\zeta=0$, $1/\lambda$ is an eigenvalue. \\

Let us fix $\lambda=\lambda_n$ and an eigenfunction $\upsilon_n(\zeta)$ such that
$$
\upsilon_n(\zeta) =(\zeta_+-\zeta)^{\frac{2\nu-1}{2}}(1+[(\zeta_+-\zeta)^2]_1)
$$
and put
\begin{equation}
w_n(z)=\kappa (\bar{\rho}\mu)^{-\frac{1}{4}} \upsilon_n(\zeta)\qquad\mbox{with}\qquad
\kappa=2^{\frac{2\nu-1}{2}}(\nu\mathsf{g}l^2)^{-\frac{\nu-1}{2}}C_{\bar{\rho}}^{\frac{1}{2}}.
\end{equation}
Then the solution $w_n$ of the equation \eqref{TG} satisfies
\begin{equation}
w_n(z)=1+[z_+-z]_1 \qquad\mbox{as}\qquad z \rightarrow z_+-0.
\end{equation}
Of course $w_n$ satisfies the usual Dirichlet boundary condition \eqref{TGBC} at the regular boundary $z=0$.

We have the associated particular solution \eqref{2.1a}\eqref{2.1b} of the Type (1) or
\eqref{Type1a}\eqref{Type1b} of the Type (2) of the linearized problem by putting
\begin{subequations}
\begin{align}
u(z)&=-\frac{1}{l}\frac{d}{dz}w_n(z), \\
w(z)&=w_n(z), \\
\delta\check{P}(z)&=-\frac{\lambda\bar{\rho}}{l^2}\frac{d}{dz}w_n(z), \\
\lambda&=\lambda_n.
\end{align}
\end{subequations}\\

Let us fix this particular  solution of the linearized problem and denote
$u=u^{L(l,n)}, w=w^{L(l,n)},\delta\check{P}=\delta\check{P}^{L(l,n)}$.
Accordingly for the Type (1) we shall denote
\begin{subequations}
\begin{align}
\mbox{\boldmath$\xi$}^{L_1(l,n)}&=
\begin{bmatrix}
u^{L(l,n)}(z)\cos lx \\
0 \\
w^{L(l,n)}(z)\sin lx
\end{bmatrix}
\sin \sqrt{\lambda}t, \\
\delta P^{L_1(l,n)}&=\delta\check{P}^{L(l,n)}(z)\sin lx\sin\sqrt{\lambda}t. 
\end{align}
\end{subequations}
and for the Type (2) we shall denote 
\begin{subequations}
\begin{align}
&\mbox{\boldmath$\xi$}^{L_2(l,,n)}=
\begin{bmatrix}
u^{L(l,n)}(z)(\cos(lx-\sqrt{\lambda}t)-\cos lx) \\
0 \\
w^{L(l,n)}(z)(\sin(lx-\sqrt{\lambda}t)-\sin lx)
\end{bmatrix}, \\
&\delta P^{L_2(l,,n)}=\delta\check{P}^{L(l,n)}(z)(\sin(lx-\sqrt{\lambda}t)-\sin lx) 
\end{align}
\end{subequations}

\subsection{}

Let us observe the movement of the vacuum boundary surface $\Gamma$, which is the boundary of the domain
$$\Pi=\{ \rho > 0\}=
\{ \bar{\rho}+\Delta{\rho} >0 \}=\{ {\bar{\rho}}+\delta{\rho}|_{t=0} >0 \}
=\{ \overset{\circ}{\rho} >0 \},
$$
described in the Eulerian co-ordinates $\underline{\mathbf{x}}=(\underline{x}^1, \underline{x}^2,
\underline{x}^3)=(\underline{x}, \underline{y}, \underline{z})$.\\

Type (1): We consider the solution $\mbox{\boldmath$\xi$}=\varepsilon\mbox{\boldmath$\xi$}^{L_1(l,n)},
\delta P=\varepsilon\delta P^{L_1(l,n)}$ with small parameter $\varepsilon, |\varepsilon|\ll 1$. Then the vacuum boundary $\Gamma$ turns out to be $\{ z=z_+\}$ since $\delta{\rho}|_{t=0}=0$, and it is 
expressed in the Eulerian co-ordinates as 
\begin{subequations}
\begin{align}
\underline{x}&=x+\varepsilon u^{L(l,n)}(z_+)\cos lx \sin \sqrt{\lambda_n}t, \label{VB1a}\\
\underline{y}&=y, \\
\underline{z}&=z_++\varepsilon w^{L(l,n)}(z_+) \sin lx \sin\sqrt{\lambda_n}t \nonumber \\
&=z_++\varepsilon \sin lx \sin\sqrt{\lambda_n}t
\end{align}
\end{subequations}
Since
\begin{align*}
\frac{\partial \underline{x}}{\partial x}(x,y,z_+,t)&=1-\varepsilon
lu^{L(l,n)}(z_+)\sin lx \sin\sqrt{\lambda_n}t \\
&=1+O(\varepsilon)
\end{align*}
uniformly, we can assume that $\partial\underline{x}/\partial x >0$ and 
\eqref{VB1a} can be solved as
\begin{equation}
x=\varphi(\underline{x}, t)=\underline{x} + O(\varepsilon),
\end{equation}
provided that $|\varepsilon|$ is sufficiently small. Then $\Gamma$ can be expressed as
\begin{align}
\underline{z}&=\underline{\mathcal{Z}}(\underline{x},\underline{y}, t)  \nonumber \\
&:=z_++\varepsilon \sin l\varphi(\underline{x}, t)\sin\sqrt{\lambda_n}t \nonumber \\
&=z_++\varepsilon(\sin l\underline{x} + O(\varepsilon))\sin\sqrt{\lambda_n}t.
\end{align}
This means that the vacuum boundary surface vibrates around the background stratospheric hight $z=z_+$ as 
$$ \underline{z}=z_++\varepsilon \sin l\underline{x}\sin\sqrt{\lambda_n}t $$
approximately modulo $O(\varepsilon^2)$.

Type (2).  We consider the solution $\mbox{\boldmath$\xi$}=\varepsilon\mbox{\boldmath$\xi$}^{L_2(l,,n)},
\delta P=\varepsilon \delta P^{L_2(l,,n)}$ with $|\varepsilon| \ll 1$. 
The vacuum boundary $\Gamma$ is the boundary of $\Pi=\{{\bar{\rho}}(\mathbf{x})+\delta{\rho}|_{t=0}(\mathbf{x})>0\}$. Since
$$\bar{\rho}(\mathbf{x})+\delta{\rho}|_{t=0}(\mathbf{x})
=\bar{\rho}(z)-\varepsilon\frac{d\bar{\rho}}{dz}w^{L(l,n)}(z)\sin lx,
$$
we see that 
$$\Pi = \{ z_+-z >0 \quad\mbox{and}\quad
1-\varepsilon\frac{1}{\bar{\rho}}\frac{d\bar{\rho}}{dz} w^{L(l,n)}(z)
\sin lx >0 \}. $$
But recall 
\begin{align*}
&-\frac{1}{\bar{\rho}}\frac{d\bar{\rho}}{dz}=\nu\frac{1}{z_+-z}(1+[z_+-z]_1), \\
&w^{L(l,n)}(z)=1+[z_+-z]_1.
\end{align*} Therefore we can find $\Upsilon(X)=1+[X]_1$ such that,
for $z_+-z>0$, 
\begin{align*}
&1-\varepsilon\frac{1}{\bar{\rho}}\frac{d\bar{\rho}}{dz} w^{L(l,n)}(z)\sin lx > 0 \\
&\Leftrightarrow \\
&z < z_++\Upsilon(\varepsilon\nu \sin lx),
\end{align*}
provided that $|\varepsilon|\nu \ll 1$. Put
\begin{equation}
\mathcal{Z}(x):=z_++\min\{\Upsilon(\varepsilon\nu\sin lx), 0\}.
\end{equation}
Then $\Pi $ turns out to be
$\{ z < \mathcal{Z}(x) \}$
so that $\Gamma$ is $\{ z=\mathcal{Z}(x) \}$.

Now the vacuum boundary $\Gamma$ can
expressed in the Eulerian co-ordinates as 
\begin{subequations}
\begin{align}
\underline{x}&=x+\varepsilon u^{L(l,n)}(\mathcal{Z}(x))(\cos(lx-\sqrt{\lambda_n}t)-
\cos lx) \label{VB4a} \\
\underline{y}&=y, \\
\underline{z}&=\mathcal{Z}(x)+ 
\nonumber \\
&+\varepsilon w^{L(l,n)}(\mathcal{Z}(x))(\sin(lx-\sqrt{\lambda_n}t)-\sin lx).
\end{align}
\end{subequations}
We see that \eqref{VB4a} can be solved monotonically as
\begin{equation}
x=\varphi(\underline{x}, t)=\underline{x}+O(\varepsilon).
\end{equation}
So $\Gamma$ can be expressed as
\begin{align}
\underline{z}&=\underline{\mathcal{Z}}(\underline{x}, \underline{y}, t)  \nonumber \\
&:=\mathcal{Z}(\varphi(\underline{x},t))+ 
\nonumber \\
&+\varepsilon w^{L(l,n)}(\mathcal{Z}(\varphi(\underline{x},t)) 
 \Big[\sin(l\varphi(\underline{x},t)-\sqrt{\lambda_n}t)-\sin l
\varphi(\underline{x}, t)\Big].
\end{align}
This means that the vacuum boundary surface is a traveling wave of the form
$$ \underline{z}=z_++(\varepsilon\nu\sin l\underline{x})\wedge 0 +\varepsilon
[\sin(l\underline{x}-\sqrt{\lambda_n}t) -\sin l\underline{x} ]
$$
approximately modulo $O(\varepsilon^2)$.

\section{Nonlinear Problem}

 We should justify the linearized approximation discussed in the preceding Section in the following sense:

{\it Let us fix the solution $\mbox{\boldmath$\xi$}^{L}, \delta P^{L}$ of the linearized problem, say,
$ \mbox{\boldmath$\xi$}^L=\mbox{\boldmath$\xi$}^{L_1(l,n)} $ or $=\mbox{\boldmath$\xi$}^{L_2(l,,n)}$, and $\delta P^L=\delta P^{L_1(l,n)}$ or $=
\delta P^{L_2(l,,n)}$ in particular. Let us suppose a finite $T$ is given. Then we should prove the existence of solutions $\mbox{\boldmath$\xi$}, \delta P$ on $0\leq t \leq T$ of the original non-linear problem of the form
\begin{equation}
\mbox{\boldmath$\xi$}=\varepsilon\mbox{\boldmath$\xi$}^L+O(\varepsilon^2), \qquad \delta P=\varepsilon\delta P^L+O(\varepsilon^2)
\end{equation}
in a suitable functional space, with $\delta\rho|_{t=0}=\varepsilon\delta\rho^L|_{t=0}$,
for any $\varepsilon$ such that $|\varepsilon|\leq \epsilon_0(T)$, $\epsilon_0(T)$
being a small positive number which may depend on $T$. }\\

Such justifications have been successfully done in \cite{TMFE}, \cite{TMOJM}, \cite{JJ}
for barotropic movements of gas with physical vacuum boundary. \\

Here we consider perturbations which are independent of the $x^2$-direction. Therefore, if
we write
\begin{equation}
x^1=x,\quad x^3=z,\quad \xi^1=X,\quad \xi^2=0, \quad \xi^3=Z,
\end{equation}
then the unknown functions $X, Z, \delta P$ are functions of $t, x, z$, and the equations
\eqref{1.15} and \eqref{NLEq} are reduced to
\begin{equation}
\mathrm{det}
\begin{bmatrix}
1+\frac{\partial X}{\partial x} & \frac{\partial X}{\partial z} \\
  &  \\
\frac{\partial Z}{\partial x} & 1+\frac{\partial Z}{\partial z}
\end{bmatrix}
=1, \label{3.3}
\end{equation}
and
\begin{equation}
\begin{bmatrix}
1+\frac{\partial X}{\partial x} & \frac{\partial X}{\partial z} \\
  &  \\
\frac{\partial Z}{\partial x} & 1+\frac{\partial Z}{\partial z}
\end{bmatrix}
\overset{\circ}{\rho}
\frac{\partial^2}{\partial t^2}
\begin{bmatrix}
X \\
 & \\
Z
\end{bmatrix}
+\mathsf{g}\delta\rho
\begin{bmatrix}
\frac{\partial X}{\partial z} \\
  &  \\
1+\frac{\partial Z}{\partial z}
\end{bmatrix}
+ 
\begin{bmatrix}
\frac{\partial}{\partial x} \\
 & \\
\frac{\partial}{\partial z}
\end{bmatrix}
\delta P =0. \label{3.4}
\end{equation}
Here we read
\begin{align}
&\delta\rho=\delta\rho(t,x,z)=-(\bar{\rho}(z+Z(t,x,z))-\bar{\rho}(z))+\delta\rho|_{t=0}(x,z), 
\label{3.5}\\
&\overset{\circ}{\rho}=\overset{\circ}{\rho}(x,z)=\bar{\rho}(z)+\delta\rho|_{t=0}(x,z). \label{3.6}
\end{align}

But we have to eliminate $\delta P$ from the equation \eqref{3.4}. This may be done by replacing \eqref{3.4} by 
\begin{align}
&\frac{\partial}{\partial z}\Big[\overset{\circ}{\rho}\Big(\Big(1+\frac{\partial X}{\partial x}\Big)
\frac{\partial^2X}{\partial t^2}+\frac{\partial X}{\partial z}\frac{\partial^2Z}{\partial t^2}\Big) +\mathsf{g}\delta\rho\frac{\partial X}{\partial z}\Big]= \nonumber \\
&=\frac{\partial}{\partial x}\Big[\overset{\circ}{\rho}
\Big(\frac{\partial Z}{\partial x}\frac{\partial^2X}{\partial t^2}+
\Big(1+\frac{\partial Z}{\partial z}\Big)\frac{\partial^2Z}{\partial t^2}\Big)+
\mathsf{g}\delta\rho\Big(1+\frac{\partial Z}{\partial z}\Big)\Big]. \label{3.7}
\end{align}

So, we should solve the system of equations \eqref{3.3} \eqref{3.7}
for unknown functions $X(t, x, z), Z(t, x, z)$ such that
\begin{equation}
X=\varepsilon(\xi^L)^1 +O(\varepsilon^2),\qquad
Z=\varepsilon(\xi^L)^3+O(\varepsilon^2)
\label{3.8}
\end{equation}
under the initial condition $\delta\rho|_{t=0}=\varepsilon\delta\rho^L|_{t=0}$.

But this task for the present system of equations is not so easy. \\

Anyway let us reflect on the linearized approximation of the system \eqref{3.3}\eqref{3.7}.

Clearly the linearized approximation  of \eqref{3.3} is
\begin{equation}
\frac{\partial X}{\partial x}+\frac{\partial Z}{\partial z}=0. \label{3.9}
\end{equation}
On the other hand, since the linearized approximations of 
\eqref{3.5}, \eqref{3.6}
are
\begin{equation}
\delta\rho=-\frac{d\bar{\rho}}{dz}Z+\delta\rho|_{t=0}, \qquad \overset{\circ}{\rho}=\bar{\rho} + \delta\rho|_{t=0},
\end{equation}
the linearized approximation of \eqref{3.7} reads
\begin{equation}
\bar{\rho}\frac{\partial^2}{\partial t^2}\Big(\frac{\partial X}{\partial z}-\frac{\partial Z}{\partial x}\Big)
+\frac{d\bar{\rho}}{dz}\frac{\partial^2 X}{\partial t^2}+
\mathsf{g}\frac{d\bar{\rho}}{dz}\frac{\partial Z}{\partial x}=0. \label{3.11}
\end{equation}

But the equation \eqref{3.9} suggests to introduce the stream function$\psi=\psi(t,x,z)$ such that
\begin{equation}
X=\frac{\partial\psi}{\partial z},\qquad Z=-\frac{\partial\psi}{\partial x}. \label{3.12}
\end{equation}
Then \eqref{3.12} reads
\begin{equation}
\bar{\rho}\frac{\partial^2}{\partial t^2}\triangle \psi+
\frac{d\bar{\rho}}{dz}\frac{\partial^2}{\partial t^2}\frac{\partial\psi}{\partial z}-
\mathsf{g}\frac{d\bar{\rho}}{dz}\frac{\partial^2\psi}{\partial z^2}=0, \label{3.13}
\end{equation}
where $\triangle$ stands for 
$\displaystyle \frac{\partial^2}{\partial x^2}+\frac{\partial^2}{\partial z^2}$.
This is the equation to be solved. In this context we see that we took
$$\psi=\frac{1}{l}w(z)\cos lx \sin \sqrt{\lambda}t $$
for Type (1) and that we took
$$\psi=\frac{1}{l}w(z)\Big(\cos(lx-\sqrt{\lambda}t)-\cos lx\Big) $$
for Type (2). Actually this taking reduces the equation \eqref{3.13} to the Taylor-Goldstein equation \eqref{TG}.


\end{document}